\documentclass[12pt]{amsart}
\usepackage{amssymb}
\usepackage{amsmath}
\usepackage[T1]{fontenc}
\usepackage[utf8]{inputenc}
\usepackage[british]{babel}
\usepackage{framed}
\usepackage{csquotes}
\usepackage[backend=biber, doi=false,isbn=false,url=false]{biblatex}
\title[Permutation groups of finite Morley rank]{Binding groups, permutation groups and modules of finite Morley rank}
\author{Alexandre Borovik}
\address{The University of Manchester, Department of Mathematics, Oxford Road, Manchester M13 9PL, United Kingdom}

\email{alexandre@borovik.net}
\author{Adrien Deloro}
\address{Sorbonne Université, Institut de Mathématiques de Jussieu--Paris Rive Gauche, \textsc{cnrs}, Université Paris Diderot.
Campus Pierre et Marie Curie, case 247, 4 place Jussieu, 75252 Paris cedex 5, France
}
\email{adrien.deloro@imj-prg.fr}

\newtheorem{conjecture}{Conjecture}
\newtheorem*{CZconjecture*}{Cherlin-Zilber Conjecture}
\newtheorem{problem}[conjecture]{Problem}
\newtheorem*{fact*}{Fact}
\newtheorem*{definition*}{Definition}
\newtheorem{theorem}{Theorem}
\theoremstyle{remark}
\newtheorem*{example*}{Example}
\newcommand{\dcl}{\operatorname{dcl}}
\newcommand{\End}{\operatorname{End}}
\newcommand{\rk}{\operatorname{rk}}

\newcommand{\PGL}{\operatorname{PGL}}
\newcommand{\PSL}{\operatorname{PSL}}
\newcommand{\GA}{\operatorname{GA}}
\newcommand{\GL}{\operatorname{GL}}
\newcommand{\SL}{\operatorname{SL}}
\newcommand{\pSL}{\operatorname{(P)SL}}
\newcommand{\Nat}{\operatorname{Nat}}

\newcommand{\bF}{\mathbb{F}}
\newcommand{\bG}{\mathbb{G}}
\newcommand{\bK}{\mathbb{K}}
\newcommand{\bN}{\mathbb{N}}
\newcommand{\bP}{\mathbb{P}}
\newcommand{\bZ}{\mathbb{Z}}

\newcommand{\ux}{\underline{x}}

\AtBeginBibliography{\small}

\begin{document}

\maketitle

{\small
\begin{quote}
\emph{If one has\/ {\rm (}or if many people have{\rm )} spent decades classifying certain objects, one is apt to forget just why one started the project in the first place. Predictions of the death of group theory in 1980 were the pronouncements of just such amnesiacs.} \hfill \cite[p.~4]{Wilson2009}\\
\end{quote}
}

\tableofcontents

The present survey aims at being a list of Conjectures and Problems in an area of model-theoretic algebra wide open for research, not a list of known results. To keep the text compact, it focuses on structures of finite Morley rank, although the same questions can be asked about other classes of objects, for example, groups definable in $\omega$-stable and $o$-minimal theories. In many cases, answers are not known even in the classical category of algebraic groups over algebraically closed fields.

\section{Introduction and background}

Groups of finite Morley rank are groups equipped with a notion of dimension which assigns to every definable set $X$ a dimension, called \emph{Morley rank} and denoted $\rk(X)$, satisfying well-known and fairly rudimentary axioms given for example in \cite{BNGroups,PStable}.
Examples are furnished by algebraic groups over algebraically closed fields, with $\rk(X)$ equal to the dimension of the Zariski closure of $X$; in particular by affine/linear algebraic groups.
However the first-order setting forbids the use of methods from algebraic geometry, be it in the modern (functorial technology), ancient (the Zariski topology on algebraic varieties), or even naive (the Jordan decomposition from linear algebra) sense of the term.

While we believe that groups of finite Morley rank form a topic legitimate \emph{per se}, more will be said on their model-theoretic relevance in \S~\ref{S:binding}.

\subsection{The classification programme}

The central conjecture on \emph{abstract} groups of finite Morley rank was formulated 25 years ago by Gregory Cherlin and Boris Zilber independently.

\begin{CZconjecture*}
Simple infinite groups of finite Morley rank are groups of points $\bG(\bK)$ of simple algebraic groups $\bG$ in algebraically closed fields $\bK$.
\end{CZconjecture*}

The theory of groups of finite Morley rank has started in pioneering works by Zilber \cite{ZGroups} and Cherlin \cite{CGroups} and was developed to a level of remarkable technical sophistication in 50+ papers by Alt\i nel, Berkman, Borovik, Cherlin, Deloro, Fr\'econ, Jaligot, Nesin.
The bulk of the work in the field is dedicated to the Cherlin-Zilber conjecture; there again, most of it is through the prism of finite group theory and the analogy with the Classification of the Finite Simple Groups, \textsc{cfsg}.
The most important result so far has been achieved in the special case leading to fields of characteristic $2$; its proof takes a 500+ page book \cite{ABCSimple}.

\begin{theorem}[\cite{ABCSimple}]\label{th:even}
If a simple group $G$ of finite Morley rank contains an infinite elementary abelian $2$-subgroup {\rm (}we say in this situation that $G$ is of \emph{at least even type}{\rm )} then it is isomorphic to a simple algebraic group over an algebraically closed field of characteristic $2$.
\end{theorem}

 In view of the persisting lack of an analogue of the Feit-Thompson ``odd order'' theorem \cite{feit1963}, the following result dramatically clarified the picture by separating reasonable (with involutions) and at the moment desperate (without involutions) situations.

\begin{theorem}[\cite{BBCInvolutions}]\label{th:degen}
If a simple group $G$ of finite Morley rank has a finite Sylow $2$-subgroup then it contains no involution.
\end{theorem}

In view of these results, the remaining configurations in a proof of the Cherlin-Zilber Conjecture involve either groups of \emph{degenerate type}, that is, simple groups without involutions (with a possibility of a counterexample emerging) or groups where a Sylow 2-subgroup contains a non-trivial divisible abelian subgroup of finite index (groups of \emph{odd type}).
A first approach to classifying groups of odd type was outlined by Borovik \cite{BSimple} and Berkman \cite{BClassical,BBIdentification} and developed in detail by Burdges in a series of papers \cite{BSignalizer,BBender}.

The present state of affairs is condensed in the following theorem which summarises a series of works by Burdges, Cherlin and Jaligot and reduces the classification of groups of odd type to a number of ``small'' configurations (further restricted in works by Deloro and Jaligot \cite{DJInvolutive}).

\begin{theorem}
Let $G$ be a %
minimal counterexample to the Cherlin-Zilber conjecture 
and assume that $G$ contains an involution.
Then %
its Pr\"{u}fer $2$-rank {\rm (}viz.~the number of copies of $\bZ_{2^\infty}$ in the direct sum decomposition $T= \bZ_{2^\infty}\times \cdots \times \bZ_{2^\infty}$ of a maximal divisible $2$-subgroup $T$ of $G${\rm )} is at most $2$.
\end{theorem}

\subsection{Concrete groups of finite Morley rank}

The analogy with finite group theory will be useful: the classification of the finite simple groups (frequently used in the form of one of its numerous corollaries, the classification of 2-transitive permutation groups) has had a profound impact on combinatorics, discrete mathematics, and representation theory.

In model theory, groups of finite Morley rank naturally appear as groups of automorphism of structures of a certain kind, and one should expect that deep structural results for groups of finite Morley rank will have a similar strong impact on model theory.

The proof of the following theorem, due to Borovik and Cherlin, is an indicator of the role of the classification technique: an answer to a basic question about actions of groups of finite Morley rank required the use of the full strength of Theorems 1 and 2 and full range of techniques developed for the study of groups of odd type.

Following \cite{BCPermutation}, we call a permutation group $(G, X)$ \emph{definably primitive} if there is no nontrivial definable $G$-invariant equivalence relation on $X$.

\begin{theorem}\label{Primitive} \emph{\cite[Theorem 1]{BCPermutation}}
There exists a function $f: \bN \rightarrow \bN$ with the following property.
If a group $G$ of finite Morley rank acts on a set $X$ transitively and definably primitively and the action is $\omega$-stable, then:
\[\rk(G) \leq f(\rk X).\]
\end{theorem}

The next example shows that the primitivity condition is 
necessary, even in the algebraic category.

\begin{example*}\label{e:matrices}
Let $\bK$ be an algebraically closed field and
$A \simeq \bK^n$ be a $\bK$-vector
space of dimension $n$. Let $T \simeq \bK^\times$
be a one-dimensional torus acting on $A$ by
matrices $\operatorname{diag}(t, t^2, \ldots, t^n)$ in a fixed basis $e_1, \ldots,
e_n$ of $A$. Let $B$ be a hyperplane in $A$
trivially intersecting the subspaces $\bK e_1, \ldots, \bK e_n$. Then
$\cap_{t \in T}{B^t} = 0$ and the right coset action of the natural
semidirect product $G=A \rtimes T$ on the coset space $X = G/B$ is faithful
and transitive. Since the stabiliser of the point $B$ is not maximal in
$G$ (because $B < A$), this action is imprimitive. However, $\dim G = n+1$ and $\dim X = 2$.
\end{example*}

Together with Theorem 4, a result by Macpherson and Pillay \cite{MPPrimitive} (who transferred to groups of finite Morley rank the classical O'Nan-Scott Theorem \cite{LPSONan} about the structure of finite permutation groups) gives hope that definably primitive groups of finite Morley rank are open to analysis.
Hence
time has come to
start to apply systematically the highly developed classification machinery of groups of finite Morley rank to questions in model theory. 

The prominent connection between the \textsc{cfsg} and the study of groups of finite Morley rank needs to be revisited at a new level. Indeed some preliminary investigations strongly suggest that a wide range of ideas from the theory of finite permutation groups and from the representation theories of finite and algebraic groups are likely to be applicable in the model-theoretic context; the reader will find more on that in subsequent sections.

\section{Binding groups and bases}\label{S:binding}

\subsection{Binding groups}

Groups of finite Morley rank emerged as binding groups in Zilber's analysis of arbitrary $\aleph_1$-categorical structures; which made them a focal point of model-theoretic algebra. By his work, any uncountably categorical structure is controlled by certain definable groups of permutations (which have finite Morley rank, by definability).
These \emph{binding groups} are model-theoretic analogues of Galois groups, and were introduced in the 1970s by Zilber \cite{Zi-binding} before being developed in other contexts by Hrushovski. 

Binding groups are very natural mathematical objects, as the following canonical example shows. Let $X$ be a finite-dimensional vector space $X$ and $G = \GL(X)$ its automorphism group.
Now fix a basis of $X$; writing in coordinates the images under $g \in G$ of the various basic vectors, that is, writing down tuples of vectors known as matrices, the equations obtained allow us to express the composition of automorphisms as products of matrices, and the group $G$ becomes \emph{definable} in $X$ (finite tuples which allow coding are a general notion in permutation group theory: a \emph{base} of a permutation group, to which we return in \S~\ref{s:base}). 
It was Zilber's seminal discovery that a similar construction can be carried out over a vast class of mathematical structures which have good logic properties of ``dependency'' between elements (although the nature of this dependency could be far away from that of linear dependency in vector spaces).

The following theorem gives a class of more general examples.

\begin{theorem}[\cite{PStable}] \label{th:binding}
Let $T$ be an $\omega$-stable theory, $M \models T$ a prime model over $\emptyset$, and let $P$, $Q$ be $\emptyset$-definable sets, with $P$ being $Q$-internal (that is, $P \subset \dcl (Q \cup F)$ for some finite $F$).
Then the group of automorphisms of $M$ which fix $Q$ pointwise induces a definable group of automorphisms of $P$, called the \emph{binding group} of $P$ over $Q$.
\end{theorem}

Even without the complete classification of groups of finite Morley rank, the results and methods already developed are powerful enough to start a systematic structural theory of binding groups, first in the context of finite Morley rank, then hopefully in the more general context of stable and simple theories.

One needs to look both at foundational works by Zilber, Poizat, Hrushovski and at a number of recent developments in the study of internality (more generally, analysability) and binding groups (not restricted to the $\omega$-stable and finite Morley rank context) in works of Ben-Yaacov, Hart, Kamensky, Moreno, Shami, Toma\v{s}i\'{c}, Wagner \cite{BTWGroup,HSType,KDefinable,KCategorical,MIterative} with the aim to identify sufficiently general sufficient conditions which would imply good combinatorial and algebraic properties of binding groups. In particular, we are interested in sufficient conditions for generic 2-transitivity of actions of binding groups on some specific orbits (see the definition in \S~\ref{sec:mt}) and, more generally, for the presence of involutions.

\subsection{Bases and parametrisation}\label{s:base}

The concept of a base of a permutation group plays a crucial role in computational methods for finite simple groups and is directly linked to binding groups in the finite Morley rank context. It is crucially important to get some control over the behaviour of bases.

\begin{definition*}
Let $(G, X)$ be a permutation group.
\begin{itemize}
\item
A subset of $X$ is said to be a \emph{base} for $G$ if its pointwise stabiliser in $G$ is trivial.
\item
The minimal cardinality of a base for $G$ is denoted by $b(G)$.
\end{itemize}
\end{definition*}
If the action of $G$ on $X$ is definable and $G$ has a descending chain condition on definable subgroups (for example, if $G$ is definable in an $\omega$-stable, or $o$-minimal structure), then it is easy to see that $b(G)$ is finite. In general a good control on the size of bases of permutation groups is an essential tool to understand parametrisations of binding groups: in the notation of Theorem~\ref{th:binding}, if $(b_1,\dots, b_n)$ is a base for the action of $G$ on $P$, then the map $G \to P^n$, $g \mapsto (b_1^g,\dots,b_n^g)$ is a definable parametrisation of $G$ suitable for analysing the action of $G$ on $P$.

Incidentally observe that this parametrisation proved to be very efficient in computations in finite permutation groups, where it was introduced by Sims \cite{Sims1971}. A detailed description can be found in Seress \cite[\S~5.2]{seress2003}. A construction of a binding group in the context of probabilitic ``black box'' computations in finite groups can be found in \cite{BY2018}.

\subsection{Sharp bounds for bases}

Let $(G, X)$ be a permutation group of finite Morley rank (i.e.~we have a faithful action, definable in some theory of finite Morley rank). Then the following holds: $\rk G \leqslant b(G)\cdot\rk X$. The matrix example above provides the sharpest bound: $\rk G = b(G)\cdot \rk X$. Moreover, in that special case $b(G) = \rk X$.

This example makes the following conjecture very natural.

\begin{conjecture}\label{conj:bases}\label{base}
In the finite Morley rank context, if $G$ is connected and definably primitive on $X$ then
\[ b(G) < c \cdot \rk X \]
for some constant $c$ which does not depend on $G$ and $X$. Moreover, the set of minimal bases is generic in $X^{b(G)}$, i.e.\ of rank $b(G)\cdot \rk X$.
\end{conjecture}
In the particular case of simple algebraic groups acting on varieties, this (and stronger statements) is a consequence of much more explicit recent results by Burness, Guralnick and Saxl. One of them is as follows.

\begin{fact*}[{\cite{burness}}]
Let $G$ be a simple algebraic group over an algebraically closed field and let $X$
be a primitive $G$-variety. Assume $G$ is not a classical group in
a subspace action. Then $b(G) \leqslant 6$.
\end{fact*}

It is instructive to compare Conjecture~\ref{conj:bases} with the following result by Liebeck which heavily depends on the \textsc{cfsg}.

\begin{fact*}[{\cite{LMinimal}}]
If $F$ is a finite primitive permutation group of degree $n$ then:
\begin{itemize}
\item[(i)]
either $b(F) < 9 \log_2 n$,
\item [(ii)]
or $F$ is a subgroup of $S_m \wr S_r$ containing $(A_m)^r$, where the action of $S_m$ is on $k$-sets and the wreath product has the product action of degree $\binom{m}{k}^r$.
\end{itemize}
\end{fact*}

But in the finite Morley rank context, for basic reasons there are no connected groups which could act in a way similar to that of clause (ii) of Liebeck's theorem. Meanwhile, clause (i) looks exactly like our Conjecture~\ref{conj:bases} if we use the usual Lang-Weil type analogy between the rank and the logarithm of cardinality of finite definable sets.

However we should not underestimate the difficulty of Conjecture~\ref{conj:bases}; most likely, all accumulated knowledge about groups of finite Morley rank will be needed for work on it.

In what follows, we outline a smaller Conjecture~\ref{c:n-trans} which, we believe, will almost inevitably arise in work on Conjecture~\ref{conj:bases}.

\section{Permutation groups of finite Morley rank}

This section and the next may be seen as having more relations with algebraic group theory than pure model theory. Let us first discuss permutation groups. Recall that a group $G$ acting on a set $X$ is \emph{$n$-transitive} if in its induced action on the set $X^{[n]}$ of $n$-tuples of distinct elements of $X$, $G$ has a unique orbit.

It is \emph{sharply} $n$-transitive if in addition, the stabiliser $G_{\ux}$ of any such tuple is trivial.
For instance, there are no infinite sharply $4$-transitive groups (proved by Tits \cite{TGroupes} and Marshall Hall Jr. \cite{HTheorem} independently, with origins in Jordan). Sharpness is not seriously discussed before \S~\ref{s:sharp}.

We shall quickly drop the ``genuine'' transitivity assumption, relaxing it into something weaker discussed in \S~\ref{sec:mt}, but let us first recall the result in the algebraic category.

\begin{fact*}[{\cite[Satz~2]{KMehrfach}}]
Multiply transitive actions of algebraic groups are known: for transitivity degree $n = 2$ and reductive $G$, equivalent to $(\PGL_{m+1}(\bK), \bP^m (F))$; for $n = 3$, equivalent to $(\PGL_2(\bK), \bP^1(\bK))$; none for $n \geq 4$.
\end{fact*}

Of course there is (at present) nothing similar for abstract groups of finite Morley rank, and we are not making it a conjecture since it is unclear whether it would be substantially easier than a large ``odd-type'' fragment of the Cherlin-Zilber conjecture.
In any case our business here is with \emph{generic} transitivity.

\subsection{Highly generically multiply transitive groups}\label{sec:mt}

\begin{definition*}
A connected group of finite Morley rank $G$ acting definably on a set $X$ of Morley degree $1$ is \emph{generically $n$-transitive} on $X$ if $G$ has a generic orbit on $X^n = X \times\cdots\times X$.
\end{definition*}

As a matter of fact the definition carries smoothly to the $\omega$-stable context, and even to the $o$-minimal case---replacing ``generic'' by ``of small codimension''.

\begin{example*}\
\begin{itemize}
\item
The general linear group $\GL_n(\bK)$ is generically $n$-transitive in its action on the vector space $\bK^n$.
\item
The affine group $\GA_n(\bK)=\bK^n \rtimes \GL_n(\bK)$ is generically $(n+1)$-transitive in its natural action on the affine space $\bK^n$.
\item
The projective group $\PGL_{n+1}(\bK)$ is generically $(n+2)$-transitive on the $n$-dimensional projective space $\bP^n(\bK)$.
\end{itemize}
\end{example*}

In the algebraic category, the classification has been obtained only recently, and only in zero characteristic \cite{PGenerically}. Notice, for example, that $E_6(F)$ has a generically $4$-transitive action, though on a large set (the quotient space of $G$ modulo a specific parabolic subgroup).

\begin{fact*}[{\cite[Theorem~1]{PGenerically}}]
Let $G$ be a simple algebraic group over an algebraically closed field of characteristic~$0$. Then the maximal degree of generic transitivity for an action of $G$ on an irreducible algebraic variety is as follows:
\[\begin{array}{c|c|c|c|c|c|c|c|c}
A_n & B_n, n\geq 3 & C_n, n\geq 2 & D_n, n\geq 4 & E_6 & E_7 & E_8 & F_4 & G_2\\\hline
n+2 & 3 & 3 & 3 & 4 & 3 & 2 & 2 & 2
\end{array}\]
\end{fact*}

In finite Morley rank one expects the same bound.

\begin{conjecture}\label{c:n-trans}
Let $G$ be a connected group of finite Morley rank acting faithfully, definably and transitively on a set $X$ of Morley rank $n$.
Then the action is at most generically $(n+2)$-transitive.
\end{conjecture}

High degree of generic transitivity is one of the principal obstacles for getting good bounds for base sizes in Conjecture \ref{base} -- this is why it is intimately linked to Conjecture \ref{c:n-trans}.

\subsubsection{Generically multiply transitive groups: the extreme case}

The following conjecture explains the reason for the bound in Conjecture~\ref{c:n-trans}.

\begin{conjecture}[{\cite[Problem~9]{BCPermutation}}] \label{c:1}
Let $G$ be a connected group of finite Morley rank acting faithfully, definably, transitively and generically $(n+2)$-transitively on a set $X$ of Morley rank $n$.
Then $(G,X)$ is definably equivalent to $(\PGL_{n+1}(\bK), \bP^n(\bK))$ for some definable field $\bK$.
\end{conjecture}

The conjecture is known for $n = 1$ \cite{HAlmost} (actually the full classification of transitive permutation groups acting on ``strongly minimal'' sets, viz.~sets of rank and degree~$1$) and more recently $n = 2$ \cite{AWRecognizing}.
(Parenthetically, Gropp \cite{GThere} had shown that the degree of sharp \emph{generic} transitivity on a set of Morley rank $2$ is at most $6$.)

The state of the art on the general case is the following result by Alt\i nel and Wiscons.

\begin{fact*}[{\cite[Theorem~A]{AWToward}}]
Let $G$ and $X$ be as in Conjecture~\ref{c:1}. Suppose that:
\begin{itemize}
\item
$G$ is $2$-transitive on $X$;
\item
there is $x \in X$ such that the permutation group $(G_x, X\setminus\{x\})$ has a definable quotient, with classes of infinite size, that has the form $(\PGL_{k+1}(\bK), \bP^k(\bK))$ for some $k$.
\end{itemize}
Then $(G, X)$ has the desired form as well.
\end{fact*}

This is at the core of an attempted inductive approach.
However analysing Conjecture~\ref{c:1} \`a la O'Nan-Scott (\cite{MPPrimitive} for groups of finite Morley rank, and \cite{MMTPermutation} for the $o$-minimal case), one should expect to run into the configuration where a stabiliser $G_\alpha$ acts not on a set, but on an abelian group, with high generic transitivity. More on this will be said in \S~\ref{s:mtmodules}.

\subsubsection{Generically multiply transitive groups: other special cases}

A classification of generically $n$-transitive groups with $n$ sufficiently large, say, $n \geq 5$, is also likely to be achievable without the complete classification of groups of odd type.
As we pointed, Popov's work remains in characteristic~$0$ and is not known in positive characteristic. But since a high degree of generic transitivity implies a large Sylow $2$-subgroup, this makes applicable a sizable arsenal of tools developed for the analysis of groups of finite Morley rank in terms of centralisers of involutions and $2$-subgroup structure, following Burdges' ``signaliser functors'' theorems \cite{BSignalizers}. So the following is not unreasonable.

\begin{conjecture}
Let $G$ be a simple group of finite Morley rank admitting a generically $4$-transitive action. Then $G \simeq \PGL_3$ or $G \simeq E_6$.
\end{conjecture}

\subsection{Sharply $n$-transitive groups}\label{s:sharp}

In this paragraph we return to \emph{genuine} $n$-transitivity, adding sharpness. As mentioned, there are no infinite sharply $k$-transitive groups for $k \geq 4$ (Jordan, Hall Jr., Tits). \emph{Finite} sharply $2$- and $3$-transitive groups were classified by Zassenhaus \cite{ZKennzeichnung}: they are essentially like $\GA_1$ and $\PGL_2$ but over \emph{near-fields}. The situation is fairly similar in the $o$-minimal case, where the near-field is either real-like, complex-like, or quaternion-like \cite{TSharplyo}.

Now in finite Morley rank and characteristic $\neq 2$, infinite near-fields are just fields, so not surprisingly infinite sharply $3$-transitive groups of finite Morley rank all have the form $\PGL_2(\bK)$ \cite{NSharply}. Notice that without model-theoretic assumptions this is trickier: Tent constructed ``exotic'' infinite sharply $3$-transitive groups \cite{TSharply}, not arising from near-fields.

The situation is even worse with sharp $2$-transitivity. Classification attempts with no additional assumption of algebraic or logical nature are hopeless \cite{RSTSharply}. Under algebraic assumptions, a sharply $2$-transitive \emph{linear} group has the desired form $\GA_1(K)$ (the near-field $K$ is not necessarily the field providing linearity) \cite{GGSharply}.
The theory of sharply $2$-transitive groups of finite Morley rank is however unfinished. The solution in so-called ``permutation characteristic $3$'' was given only recently \cite{ABWSharply}. The relevant conjectures are stated there, if not earlier in the literature.

Since any sharply $2$-transitive group is a Frobenius group, this also explains why the classification of Frobenius groups of finite Morley rank is still open.

\section{Modules of finite Morley rank}

As opposed to the previous section which dealt with permutation groups in first-order logic, the focus is now more on actions on abelian groups (as a consequence of the O'Nan-Scott analysis, but also as interesting \emph{per se}), a topic we call \emph{first-order representation theory}. Following general model-theoretic analogy, we do not represent our algebraic structures in vector spaces but in abelian groups---called modules---hoping that the model theory will as often induce definable coordinatisation results. In a sense, Zilber's celebrated ``Field Theorem'' may be seen as the starting point of first-order representation theory.
Although one could represent other algebraic structures (with strong interest in associative rings and Lie rings, which suggests in due time moving from the definable to at least the $\vee$-definable category since enveloping structures need no longer be properly definable), we shall in the present survey be content with definable groups.

A \emph{module of finite Morley rank} consists in a group $G$ and a $G$-module $V$, all relevant structure being definable in some theory of finite Morley rank.

\subsection{The torsion-free case} \label{sec:torsion-free}

The case where $V$ is torsion-free is sufficiently well understood due to two fundamental results. The first enables linearisation \emph{a priori} and is undergoing current generalisation way beyond the finite Morley rank category \cite{DWLinearisation}. (There is something similar in positive characteristic only if one assumes $C_{\End(V)}(G)$ to be infinite.)

\begin{fact*}[{\cite{LWCanada}}]
Let $(G, V)$ be a faithful, irreducible module of finite Morley rank where $G$ is infinite and $V$ is torsion-free. Then there is a definable field over which $V$ is a finite-dimensional vector space and $G \hookrightarrow \GL(V)$, definably.
\end{fact*}

The second provides a strong division between the algebraic world and (potential) highly pathological configurations. (The situation is good although quite non-trivial in positive characteristic: see \cite{PQuelques}.)

\begin{fact*}[\cite{MPPrimitive, PQuelques}; for involutions \cite{BBDefinably, DWGeometry}]
Let $G$ be an infinite, definably linear, simple group of finite Morley rank, represented as
a definable subgroup of $\GL_n(\bK)$
for some field $\bK$ of finite Morley rank \emph{of characteristic zero}.
Then one of the following occurs.
\begin{itemize}
\item $G$ is Zariski-closed, hence a linear algebraic group.
\item
$G$ has no involutions; its Borel subgroups are conjugate self-normalising abelian subgroups of $G$; and the same applies to any definable, connected subgroup of $G$.
\end{itemize}
Furthermore, $G$ is the union of its Borel subgroups.
\end{fact*}

As a consequence, a simple algebraic group acting on a finite Morley rank module in characteristic $0$ acts algebraically; see \S~\ref{s:Steinberg} for further discussion.

\subsection{Simultaneous identification}\label{s:simultaneous}

As opposed to the next subsection, the present deals with the case where neither $G$ nor $V$ is known, and we seek simultaneous identification under assumptions on the action. There are several natural types of assumptions, but bear in mind that as opposed to the characteristic~$0$ case (\S~\ref{sec:torsion-free}), we have no a priori linearisation results.

\subsubsection{A basic case: transitive modules}

\begin{conjecture}\label{c:afterKnop}
Let $(G, V)$ be a module of finite Morley rank where $G$ is connected, and simple modulo its at most finite centre. Suppose that $G$ is transitive on $V\setminus\{0\}$ and generically $2$-transitive. Then $(G, V)$ is definably equivalent to $(\SL_n(\bK), \bK^n)$ for some definable field $\bK$.
\end{conjecture}

This is known \emph{in the algebraic category} \cite[Satz~1]{KMehrfach}, where \emph{no abelianity of $V$ is required}. Conjecture~\ref{c:afterKnop} is a basic problem where one could test new linearisation methods.

\subsubsection{Highly generically multiply transitive modules}\label{s:mtmodules}

Running Conjecture~\ref{c:1} in the category of modules of finite Morley rank gives the following.

\begin{conjecture}[{\cite[Problem~13]{BCPermutation}}]\label{c:2}
Let $(G, V)$ be a faithful module of finite Morley rank where $V$ has rank $n$ and $G$ is generically $n$-transitive.
Then $(G, V)$ is definably equivalent to $(\GL_n(\bK), \bK^n)$ for some definable field $\bK$.
\end{conjecture}

\paragraph{\it Supporting evidence.}

The conjecture is established for small values of $n$, and under sharpness assumptions.

\begin{fact*}[Deloro \cite{DActions}]\label{fact:Deloro-rank-2}
Let $(G, V)$ be a faithful, irreducible module of finite Morley rank where $G$ is connected and $V$ has rank $2$.
Then $G$ is one of the groups $\SL_2(\bK)$ and $\GL_2(\bK)$ in their natural action on $V= \bK^2$.
\end{fact*}

\begin{fact*}[Borovik and Deloro \cite{BDRank}, using Fr\'{e}con \cite{FSimple}]\label{f:BD}
Let $(G, V)$ be a faithful, irreducible module of finite Morley rank where $G$ is connected, non-soluble and $V$ has rank $3$.
Then:
\begin{itemize}
\item[(a)]
either $G = \PSL_2(\bK) \times Z(G)$ where $\PSL_2(\bK)$ acts in its adjoint action on $V\simeq \bK_+^3$,
\item[(b)]
or $G = \SL_3(\bK) \ast Z(G)$ in its natural action on $V \simeq \bK_+^3$.
\end{itemize}
\end{fact*}

Interestingly enough, the proof of the latter involves ideas from more or less all directions explored over almost forty years of the theory of groups of finite Morley rank: the bulk of the classification programme of course (both in even and odd characteristics), but also more ``hardcore'' model-theoretic tools touched upon in \S~\ref{sec:torsion-free} and relating to linearisation \cite{LWCanada} or definably linear groups \cite{PQuelques}.

Conjecture~\ref{c:2} is also confirmed in the important special case when the action of $H$ on $V$ is \emph{sharply} generically transitive \cite{BBGenerically}. There is hope that the general case could be reduced to this special result if one obtains good bounds for the Morley ranks of irreducible modules for some specific finite groups, see Problem~\ref{problem-finite-groups} and its discussion in \S~\ref{sec:modules-finite-groups}.

\paragraph{\it Classification of pseudoreflection groups.}\label{sec:pseudoreflection}

Finally, Conjecture~\ref{c:2} is likely to follow from yet another Conjecture; the latter deals with a configuration which played a prominent role in the classification of groups of even type \cite{ABCSimple}.%

\begin{conjecture}[{\cite[Problem~17]{BCPermutation}}]\label{c:pseudoreflection}
Let $(G, V)$ be a faithful, irreducible module of finite Morley rank where $G$ is connected.
Assume that $G$ contains a \emph{pseudoreflection subgroup}, that is, an abelian subgroup $R$ such that $V = [V,R] \oplus C_V(R)$ and $R$ acts transitively on the set of non-trivial elements in $[V,R]$.
Then there is a definable field $\bK$ such that $(R, [V, R)]$ is definably equivalent to $(\bK^\times, \bK_+)$, and $(G, V)$ to $(\GL(V), \bK^n)$.
\end{conjecture}

Notice that the most important case is $\rk [V,R]=1$ (and this is the only case needed for proving Conjecture~\ref{c:2}). It is also the most interesting to a finite group theorist's taste: this is the point where the \emph{classical involutions} theme \cite{ACharacterization,BClassical} comes into the plot.
Significant progress has already been achieved in \cite{BBPseudoreflection}.

\subsubsection{Relations with the classification programme}

Of course the interplay between abstract and concrete groups is subtle and goes both ways. We suggest two reasons why knowledge on finite Morley rank modules may influence the classification programme.%

\paragraph{\it The Cherlin-Zilber Conjecture for groups admitting a definable module.}

The following configuration arises at early stages of the analysis in Conjecture~\ref{base} (it appears inside primitive permutations groups of affine type \cite{MPPrimitive}) and is one of the first to be treated.

\begin{conjecture}\label{c:abelian-p-group}
Let $(G, V)$ be a faithful module of finite Morley rank, where $G$ is simple and $V$ has exponent $p > 2$.
Then $G$ is either of degenerate type, or isomorphic to an algebraic group over an algebraically closed field of characteristic $p$.
\end{conjecture}

The entire programme for the study of groups of odd type could be run in this special case. Wagner's theorem on bad fields \cite{WFields} provides a supply of $q$-tori (that is, divisible abelian $q$-groups) for primes $q\ne p$, while the seminal result by Cherlin \cite{CGood} on conjugacy of maximal ``good tori'' and subsequent work by Burdges and Cherlin \cite{BCSemisimple} will give us an efficient control over generic elements in $G$.

Notice that there is little hope to actually produce such a module out of pure model-theory: to our knowledge, there is no way to introduce a Lie algebra \emph{except presumably in the case of Zariski geometries} \cite{Hrushovki-Zilber1996,zilber_2010}.

\begin{conjecture}
Let $G$ be a simple group of finite Morley rank definable in a Zariski geometry. Then $G$ has a Lie algebra, and is a simple algebraic group.
\end{conjecture}

\paragraph{\it Quadratic pairs.}

\emph{Quadratic pairs} were introduced by Thompson in unpublished work \cite{TQuadratic}. Classically a quadratic pair is a pair $(G, M)$ where $G$ is a group, $M$ a faithful, irreducible $G$-module, and $G$ is generated by its quadratic elements, viz.~those satisfying $(g-1)^2 = 0$ in $\End(M)$.

Thompson \cite{TQuadratic} identified the \emph{finite} groups $G$ belonging to a quadratic pair of characteristic $\geq 5$: such a pair has the form $$(G_1 \ast \dots \ast G_d, M_1 \otimes \dots \otimes M_d)$$ where each $(G_i, M_i)$ is a quadratic pair with $G_i$ a quasi-simple Chevalley group.
This result is of essential use in the third generation of the \textsc{cfsg} (based on the ``amalgam method'') and one should hope that a similar classification could be useful toward the Cherlin-Zilber conjecture.

(Parenthetically the characteristic $3$ case, tackled by Ho \cite{HQuadratic} and completed by Chermak \cite{CQuadratic} \emph{using the \textsc{cfsg}}, is far more delicate. Also notice that Thompson classified the \emph{groups} $G$, not the \emph{pairs} $(G, M)$: explicit determination of quadratic modules for finite groups of Lie type was later achieved in \cite{PSQuadratic}.)

\begin{problem}
Classify quadratic pairs of finite Morley rank.
\end{problem}

\subsection{Model-theoretic representations}

From now on we suppose that $G$ is already known, for instance because it is a finite group, or an algebraic group.

\subsubsection{Model-theoretic representations of finite groups} \label{sec:modules-finite-groups}

\begin{problem} \label{problem-finite-groups}
Let $(G, V)$ be a faithful, irreducible module of finite Morley rank where $G$ is finite.
Find a lower bound on $\rk V$ in terms of the complexity of $G$.
\end{problem}

So far the cases of the symmetric and alternating groups have been determined in \cite{CDWMinimal} (reviving classical work of Dickson \cite{DRepresentations} and following \cite{BCPermutation}).
Obtaining good bounds on model-theoretic representations of extraspecial $p$-groups should give some control on representations of ``large'' algebraic groups (where the Weyl group is rich), and provide key results towards Conjecture~\ref{c:pseudoreflection}.

It is actually expected that the context of finite-dimensional theories as expounded in \cite{DWLinearisation} should be enough to develop this topic.

\subsubsection{Model-theoretic representations of algebraic groups}\label{s:Steinberg}

One could also start with known $G = \bG(\bK)$ and try to determine its modules of finite Morley rank with the hope that it is only a reasonable extension of the rational category. The following is modelled after Steinberg's ``tensor product'' theorem \cite{SRepresentations}. Notice the presence of \emph{definable} field twists.

\begin{conjecture}[{\cite{DRegard}}]\label{c:Steinberg}
Let $(G, V)$ be a faithful, irreducible module of finite Morley rank where $G = \bG(\bK)$ is the group of $\bK$-points of a simple algebraic group. Then there are:
\begin{itemize}
\item
on $V$ a compatible, definable $\bK$-vector space structure,
\item
irreducible \emph{algebraic} modules $W_1, \dots, W_d$,
\item
$(G\ltimes V)$-definable field automorphisms $\varphi_1, \dots, \varphi_d$,
\end{itemize}
such that $V \simeq \bigotimes_{i = 1}^d W_i^{\varphi_i}$ as $\bK[G]$-modules. (We say that $V$ is in the \emph{definable twist-and-tensor category}.)
\end{conjecture}

As a consequence of the results in \S~\ref{sec:torsion-free}, this is known in characteristic~$0$ \cite[Lemma~1.4]{CDSmall}, where there are no field twists, and only one term in the tensor product: the module is an algebraic representation. (Something similar holds in the $o$-minimal world \cite[Proposition~4.1]{MMTPermutation}.)
In finite Morley rank and positive characteristic, the most advanced result to-date is the following.

\begin{theorem}[\cite{DStudy}, unpublished]
If $G = \pSL_2(\bK)$ and $\rk V \leq 4 \rk \bK$, then Conjecture~\ref{c:Steinberg} holds.
\end{theorem}

So far only the case $\rk V \leq 3 \rk \bK$ is published \cite{DActions, CDSmall}; these however only deal with the natural and adjoint modules, where no twists occur. The unpublished \cite{DStudy} tackles $\Nat \SL_2(\bK) \otimes \Nat \SL_2(\bK)^\varphi$ at a considerable computational cost. As a matter of fact one could hope to push the method to $\rk V \leq 5 \rk \bK$, but for serious geometric obstructions explained in \cite{DSymmetric2}, there is at present no general strategy even for $G = \pSL_2(\bK)$. The first author however recommends studying $V$ as a $\bG(\bF_q)$-module for increasing values of $q$, which naturally brings us back to the topic of representations of finite Morley rank of finite groups as in \S~\ref{sec:modules-finite-groups}. Other hopes are formulated in the conclusion of \cite[\S~4]{DRegard}.

\section*{Acknowledgements}

The first author had many useful discussions with Gregory Cherlin, Omar Leon Sanchez and Dugald Macpherson, and wishes to thank them. 

Parts of this paper were written when both authors enjoyed the warm hospitality offered to them at the Nesin Mathematics Village in \c{S}irince, Turkey, a proud recipient of the Leelavati Prize of the International Mathematics Union.  

\printbibliography

\end{document}